\documentclass[10pt]{amsart}
\usepackage{amsmath,amscd}
\usepackage{amsmath,amsfonts}
\usepackage[mathscr]{eucal}
\usepackage{amssymb,amsmath}
\usepackage{latexsym}
\usepackage{amsfonts}
\usepackage{amscd}
\usepackage{amsthm}
\usepackage{graphics}
\usepackage{epsfig}
\usepackage{esint}
\usepackage{pdfsync}
\usepackage{xcolor}
\usepackage[pagebackref]{hyperref}
\hypersetup{
    colorlinks = true,
    linkbordercolor = {green},    
}

\renewcommand{\d}{\mathrm{d}}

\newcommand{\bbE}{{\mathbb E}}
\newcommand{\bbR}{{\mathbb R}}
\newcommand{\bbC}{{\mathbb C}}

\newcommand{\E}{{\mathrm e}}
\newcommand{\iC}{{\mathrm i}}

\newcommand{\SO}{\operatorname{SO}}

\newcommand{\Or}{\operatorname{O}}

\newcommand{\la}{\langle}
\newcommand{\ra}{\rangle}

\newcommand{\w}{{\mathchoice{\,{\scriptstyle\wedge}\,}{{\scriptstyle\wedge}}
      {{\scriptscriptstyle\wedge}}{{\scriptscriptstyle\wedge}}}}

\newcommand{\be}{\begin{equation}}
\newcommand{\ee}{\end{equation}}
\newcommand{\bpm}{\begin{pmatrix}}
\newcommand{\epm}{\end{pmatrix}}

\numberwithin{equation}{section}

\newtheorem{theorem}{Theorem}

\theoremstyle{remark}
\newtheorem{definition}{Definition}

\newtheorem{remark}{Remark}

\begin{document}

\author[R. Bryant]{Robert L. Bryant}
\address{Duke University Mathematics Department\\
         P.O. Box 90320\\
         Durham, NC 27708-0320}
\email{\href{mailto:bryant@math.duke.edu}{bryant@math.duke.edu}}
\urladdr{\href{http://www.math.duke.edu/~bryant}%
         {http://www.math.duke.edu/\lower3pt\hbox{\symbol{'176}}bryant}}

\title[Conformal Volume of Tori]
      {On the conformal volume\\
            of $2$-tori}

\date{March 17, 2025}

\begin{abstract}
This note (originally from 2015) provides 
a proof of a 1985 conjecture of Montiel and Ros 
concerning the conformal volume of tori.

This updated version adds a proof of the claim made 
in Remark~\ref{rem: confvol} below 
about the value of the conformal volume of tori in the cases
not covered by the conjecture of Montiel and Ros.
\end{abstract}



\thanks{
Thanks to Duke University for its support via a research grant 
and to the National Science Foundation 
for its support via the grant DMS-8352009 
(during which, in 1984--5, nearly all of this work was done) 
and the grant DMS-1359583 (originally DMS-1105868), 
during which this manuscript was completed.   
\hfill\break
\hspace*{\parindent} 
This is Draft Version~$0.2$.
}

\maketitle

\setcounter{tocdepth}{2}
\tableofcontents

\section{Introduction}\label{sec: intro}

\subsection{Conformal volume}\label{ssec: confvol}
In a celebrated paper~\cite{Li-Yau}, 
Li and Yau introduced and studied the notion of \emph{conformal volume} 
of a compact conformal manifold~$\bigl(M,[g]\bigr)$.

To recall this definition, let~$G_n$ be the group
of conformal automorphisms of the standard unit $n$-sphere~$\bigl(S^n,g_0\bigr)$;
it is known to be isomorphic to $\Or(n{+}1,1)$. 
Given a smooth mapping~$f:M^m\to S^n$, let $V(f)$ 
denote the volume of the (possibly degenerate) metric~$f^*(g_0)$ on~$M$.
One then defines the \emph{conformal volume} of~$f$ to be 
$$
V_c(f) = \sup_{\psi\in G_n} V(\psi{\circ}f).
$$
(This is also what Gromov~\cite{Gromov} 
has called the \emph{visual volume} of~$f$.)
If $f$ is an immersion at some point, then $V_c(f)\ge V(S^m,g_0)>0$.
Indeed, if $m>1$ (or $m\ge 1$ if $M$ is compact), 
one has~$V_c(f) > V(S^m,g_0)$ 
unless $f(M)$ lies in a totally umbilic $m$-sphere in $S^n$. 
See \cite[Proposition 1]{Bryant}.)

Now, every conformal automorphism $\psi:S^n\to S^n$ 
extends uniquely to a conformal automorphism 
of the $(n{+}1)$-ball~$B^{n+1}\subset\bbE^{n+1}$
of which $S^n$ is the boundary. 
If $\d V_f$ is the volume form on $M$ 
induced by~$f:M^m\to S^n$, then one finds that
$$
\d V_{\psi{\circ}f} 
= \frac{\bigl(1-\psi(0){\cdot}\psi(0)\bigr)^{m/2}}
        {(1-\psi(0){\cdot}f)^{m}} \d V_f\,.
$$ 
Thus, one has the alternative formula
$$
V_c(f) = \sup_{a\in B^{n+1}} \ \ 
\int_M \frac{(1-a{\cdot}a)^{m/2}}{(1-a{\cdot}f)^{m}}\ \d V_f\,.
$$

When $M^m$ is endowed with a conformal structure~$[g]$, 
Li and Yau define the \emph{conformal $n$-volume} of~$\bigl(M,[g]\bigr)$ to be
$$
V_c\bigl(M,[g],n\bigr) = \inf\left\{\ V_c(f)\ \mid f:\bigl(M,[g]\bigr)
\to \bigl(S^n,[g_0]\bigr)\ \text{is branched conformal}\ \right\}.
$$
Of course, $V_c\bigl(M,[g],n\bigr)$ is only defined when $n$ is large enough 
for there to exist a branched conformal immersion $f:M^m\to S^n$.  
Since one has~$V_c\bigl(M,[g],n\bigr)\ge V_c\bigl(M,[g],n{+}1\bigr)$
when both quantities are defined, 
one can define the \emph{conformal volume} of~$\bigl(M,[g]\bigr)$ to be
$$
V_c\bigl(M,[g]\bigr) = \lim_{n\to\infty} V_c\bigl(M,[g],n\bigr).
$$
Naturally, in the compact case, one has $V_c\bigl(M^m,[g]\bigr)\ge V(S^m,g_0)$.

In~\cite{Li-Yau}, 
Li and Yau give a number of applications of this notion of conformal volume,
showing, in particular, how it is related to the Willmore problem, 
eigenvalue estimates, 
and a number of other interesting differential geometric issues.  
They developed methods of estimating the conformal volume,
but were able explicitly to compute it only in a small number of cases.

\subsection{The case of 2-tori}\label{ssec: confvol2tori}
Montiel and Ros~\cite{Montiel-Ros} carried out a careful study 
of the conformal volume in the case of a torus~$T_\tau = \bbC/\Lambda_\tau$ 
where~$\Lambda_\tau$ is the discrete lattice in $\bbC\simeq\bbR^2$ 
generated by~$1$ and~$\tau = x + \iC y$, 
where $0\le x\le \tfrac12$ and $y\ge \sqrt{1-x^2}$,
and the metric on the torus 
is~$|\d z|^2 = \d z{\circ}\d \bar z = \d u^2 + \d v^2$.
(Every metric on the $2$-dimensional torus 
is conformal to this metric on~$\bbC/\Lambda_\tau$ 
for some unique~$\tau = x + \iC y$ satisfying these inequalities.)  

For this conformal torus, 
they constructed a conformal (in fact, homothetic) 
embedding~$f_\tau:T_\tau\to S^5$ 
with volume 
\be\label{eq: MRtorusvolume}
V(f_\tau) = \frac{4\pi^{2}y}{y^2+x^{2}-x+1}\,,
\ee
and they showed that 
\be\label{eq: MRineq}
V_c\bigl(T_\tau,\bigl[|\d z|^2\bigr]\bigr)
 \ge \frac{4\pi^{2}y}{y^2+x^{2}-x+1}\,.
\ee  

Since, the right hand side of \eqref{eq: MRtorusvolume} 
tends to zero as $y$ goes to $\infty$, 
the inequality~\eqref{eq: MRineq} 
must become strict for~$y$ sufficiently large.  

Montiel and Ros showed~\cite[Theorem 8]{Montiel-Ros} that, 
when $x=0$ and $1\le y\le \sqrt{5/3}$, equality
does, in fact, hold in~\eqref{eq: MRineq}.  
Moreover, they conjectured that equality holds in~\eqref{eq: MRineq} 
as long as $(x{-}\tfrac12)^2+y^2<\tfrac94$, i.e., $|\tau{-}\tfrac12|<\tfrac32$.
(See Remark~\ref{rem: MRconj} for their motivation.)

In this note (see Theorem~\ref{thm: MRconj}), 
I prove (a slightly stronger version of) their conjecture, i.e.,  
\be\label{eq: MReq}
V_c\bigl(T_\tau,\bigl[|\d z|^2\bigr]\bigr)= 
\frac{4\pi^{2}y}{y^2+x^{2}-x+1}
\qquad\text{when}\ |\tau{-}\tfrac12|\le\tfrac32.
\ee

The proof proceeds in two steps: 
First, I show, as did Montiel and Ros, 
that, for all $\tau$, 
one has $V_c\bigl(T_\tau,\bigl[|\d z|^2\bigr]\bigr)\ge V(f_\tau)$.  
Second, I show that $V_c(f_\tau) = V(f_\tau)$ 
when $|\tau{-}\tfrac12|\le\tfrac32$, i.e., that, for such~$\tau$,
\be\label{eq: RLBbd}
\int_{T_\tau} \frac{(1-a{\cdot}a)}{(1-a{\cdot}f_\tau)^{2}}\ \d V_{f_\tau} 
\le \frac{4\pi^{2}y}{y^2+x^{2}-x+1}
\qquad\text{for all $a\in \bbE^6$ with $a{\cdot}a < 1$.}
\ee
Since, by definition, 
$$
V_c(f_\tau)\ge 
V_c\bigl(T_\tau,\bigl[|\d z|^2\bigr],n\bigr)\ge V_c\bigl(T_\tau,\bigl[|\d z|^2\bigr]\bigr)
$$ 
for $n\ge 5$, the conjecture follows.
 
\begin{remark}[Motivation for the conjecture]\label{rem: MRconj}
By a Taylor series expansion, Montiel and Ros 
showed that the function of~$a$ on the left hand side of \eqref{eq: RLBbd}
attains a local maximum at $a=0\in\bbE^6$ when $|\tau{-}\tfrac12|<\tfrac32$ 
and that it does not have a local maximum there when $|\tau{-}\tfrac12|>\tfrac32$.

However, they were not able to estimate the integral accurately
enough to show that $a=0$ actually is the maximum of the function 
on the unit ball in $\bbE^6$ when $|\tau{-}\tfrac12|<\tfrac32$.  
My contribution to proving their conjecture is to show that, 
when $|\tau{-}\tfrac12|\le\tfrac32$, the function on the left
hand side of \eqref{eq: RLBbd} does indeed attain its supremum
on the unit ball in~$\bbE^6$ at $a=0$.
\end{remark}

\section{Balanced Volume}\label{sec: balvol}
Throughout this section, $\tau\in\bbC$ denotes $x+\iC\,y$
where $0\le x\le \tfrac12$ and $y\ge \sqrt{1-x^2}$,
and $\Lambda\subset\bbC$ is the lattice generated by $1$ and $\tau$.
Let $T=\bbC/\Lambda$ and let $\d A = \tfrac\iC2\,\d z\w\d\bar z$
be the area form on~$T$.  (I will write $z = u+\iC v$ 
to denote the real and imaginary parts of the coordinate $z$ on $\bbC$.  
While $z$ is not a well-defined coordinate on $T$, 
the differential $\d z$ is well-defined on $T$, and, consequently, 
for any vector-valued function~$f$ on~$T$, 
one has $\d f = f_{z}\,\d z + f_{\bar z}\,\d \bar z$, so that
the quantities $f_{z}$ and $f_{\bar z}$ are also well-defined.)

A $C^\infty$ map $\phi:T\to S^n$ is said to be \emph{$\d A$-balanced} if
\be
\int_{T} \phi\>\d A = 0\in\bbE^{n+1}.
\ee
The map $\phi$ is said to be \emph{weakly conformal} 
if $\phi_{z}\cdot\phi_{z} = 0$.  When $\phi$ is weakly conformal, 
set
\be
V(\phi) = 2\int_{T} \phi_{z}{\cdot}\phi_{\bar z}\ \d A.
\ee
Note that, when $\phi$ is weakly conformal, one has $V(\phi)>0$ 
unless $\phi$ is constant.

\begin{definition}
The \emph{balanced volume} of $T$ is the infimum $V_{b}(n,T)$
of the numbers $V(\phi)$ where $\phi$ ranges over the set of 
$\d A$-balanced, weakly conformal smooth maps $\phi:T\to S^{n}$.
\end{definition}

The following result was proved by Montiel and Ros 
\cite[Proposition 5]{Montiel-Ros} (and, independently, by me) in 1985:

\begin{theorem}\label{thm: Vbcomputations}
The following formulae hold
\begin{itemize}
\item $V_{b}(2,T) = 8\pi$ for all $\tau = x+\iC\,y$.
\item $V_{b}(n,T) = {\displaystyle\frac{4\pi^{2}y}{y^2+1}}$ 
                    for all $n\ge 3$ if $x=0$.
\item $V_{b}(n,T) = {\displaystyle\frac{4\pi^{2}y}{y^2+x^{2}-x+1}}$ 
                    for all $n\ge5$.
\end{itemize}
Moreover, there is a $\d A$-balanced, weakly conformal map
$f_\tau:T\to S^{5}\subset S^{n}$, unique up to rigid motions in $S^{n}$, 
for which $V(f_\tau)=V_{b}(5,T)$.  This $f_\tau$ is homothetic and,
when $x=0$, it is a linearly full immersion in an $S^{3}\subset S^{5}$ 
while, when $x>0$, it is a linearly full immersion in~$S^{5}$.
\end{theorem}

\begin{remark}
I do not know the values of $V_{b}(3,T)$ or $V_{b}(4,T)$ when $x>0$,
but, since $V_{b}(n,T)\ge V_{b}(n{+}1,T)$, one has the inequalities
$$
V_{b}(3,T)\ge V_{b}(4,T)\ge V_{b}(5,T) = \frac{4\pi^{2}y}{y^2+x^{2}-x+1}.
$$
\end{remark}

\begin{proof}
The case $n=2$ is special 
because a $\d A$-balanced, weakly conformal map $\phi:T\to S^{2}$ 
is either holomorphic or anti-holomorphic and hence $V(\phi) = 4|d|\pi$ 
where $d$ is the degree of~$\phi$, which necessarily satisfies $|d|\ge2$. 
(The degree $d$ cannot be zero, since, in this case,
$\phi$ is constant and hence cannot be $\d A$-balanced.  
Moreover, $d\not=\pm1$ since such a map would necessarily 
be a biholomorphism between $T$ and $S^2$.)  
Thus, $V(\phi)\ge 8\pi$, so $V_{b}(2,T) \ge 8\pi$.  
On the other hand, every torus $T$ 
does admit a holomorphic mapping~$\phi:T\to S^{2}$ with degree $d=2$.  
Moreover, it is possible to find such a $\phi$ that is~$\d A$-balanced.  
Thus, $V_{b}(2,T) = 8\pi$.

Returning to general~$n$, let $\phi:T\to S^n$ 
be a $\d A$-balanced, weakly conformal mapping.
The proof will now proceed by expanding $\phi$ in its Fourier series.

To do this, let $\Lambda^{*}\subset\bbC$ 
be the lattice dual to $\Lambda$ via the pairing
$$
\la\xi,z\ra = \Re(\xi\bar z) = \tfrac12(\xi\bar z + \bar\xi z).
$$
Then $\Lambda^{*}$ is generated as a lattice by $\{\iC/y, 1-\iC\,x/y\}$.

The Fourier series expansion of $\phi$ then takes the form
$$
\phi = \sum_{\xi\in\Lambda^{*}} \phi_{\xi}\,\E^{2\pi\iC\la\xi,z\ra}
$$
where $\phi_{0}=0$ (since $\phi$ is $\d A$-balanced) 
and $\phi_{-\xi}=\overline{\phi_{\xi}}$ (since $\phi$ is real-valued).  
Since $\phi\cdot\phi\equiv1$, one has the identities
$$
1 = \sum_{\xi\in\Lambda^{*}} |\phi_{\xi}|^{2}
\qquad\text{and}\qquad
0 = \sum_{\xi\in\Lambda^{*}} \phi_{\xi}\cdot\overline{\phi_{\xi-\eta}}
$$
for all $\eta\in\Lambda^{*}\setminus\{0\}$.
The equation $\phi_{z}\cdot\phi_{z}\equiv0$ (weak conformality) becomes
$$
0 = \sum_{\xi\in\Lambda^{*}} 
       \xi(\xi{-}\eta)\,\phi_{\xi}\cdot\overline{\phi_{\xi-\eta}}
$$
for all $\eta\in\Lambda^{*}\setminus\{0\}$.  
In particular, setting $\eta=0$, one finds%
\footnote{It will be seen that $V(\phi)\ge V_{b}(5,T)$ 
holds for all $\d A$-balanced maps~$\phi:T\to S^{n}$ 
that satisfy \eqref{eq: weak conformal average}, 
which, of course, is equivalent to the very weak condition 
 $\int_{T}(\phi_{z}{\cdot}\phi_{z})\,\d A = 0$.}
\be\label{eq: weak conformal average}
0 = \sum_{\xi\in\Lambda^{*}} 
       \xi^{2}\,|\phi_{\xi}|^{2}.
\ee

Meanwhile, computation yields
$$
V(\phi) = 2\pi^{2}y\,\sum_{\xi\in\Lambda^{*}} 
                            |\xi|^{2}\,|\phi_{\xi}|^{2},
$$
and it is this latter quantity for which I seek to determine a lower bound.

Let 
$$
\Lambda^{*}_{+} = \left\{ \xi\in\Lambda^{*}\ \vrule\ \Re(\xi)>0
\ \text{or} \left( \Re(\xi)=0\ \text{and}\ \Im(\xi)>0\right)\right\},
$$
so that $\Lambda^{*} = \Lambda^{*}_{+}\cup\{0\}\cup(-\Lambda^{*}_{+})$.

For $\xi \in \Lambda^{*}_{+}$, set $p_{\xi} = 2|\phi_{\xi}|^{2}$.
Thus, one seeks the minimum value of 
$$
L = \sum_{\xi\in\Lambda^{*}_{+}} |\xi|^{2}\,p_{\xi}
$$
subject to the conditions that $p_{\xi}\ge0$ when $\xi\in\Lambda^{*}_{+}$
while
$$
1 = \sum_{\xi\in\Lambda^{*}_{+}} p_{\xi}
\qquad\text{and}\qquad
0 = \sum_{\xi\in\Lambda^{*}_{+}} \xi^{2}\,p_{\xi}.
$$

Geometrically, this amounts to the following:  
For $\xi\in\Lambda^{*}_{+}$, define
$$
\hat \xi = \bigl(\Re(\xi^{2}), \Im(\xi^{2}), |\xi|^{2}\bigr)
\in N_+\subset \bbR^{3},
$$
where $N_+ = \left\{ \bigl(a,b,\sqrt{a^2{+}b^2}\bigr)\in\bbR^3
 \mid\ (a,b)\not=(0,0)\ \right\}$ 
is the upper right circular cone in~$\bbR^3$, 
and let $C$ be the convex hull in $\bbR^{3}$ of the set $\left\{\hat\xi\ \vrule
\ \xi\in\Lambda^{*}_{+}\right\}\subset N_+$.  Thus,
$$
C = \left\{\ \sum_{\xi\in\Lambda^{*}_{+}} p_{\xi}\,\hat\xi\ \ \vrule
\ \ p_{\xi}\ge0
\ \text{and}\ \sum_{\xi\in\Lambda^{*}_{+}} p_{\xi} = 1\ \right\}.
$$
Then one seeks the infimum of the numbers~$\ell\in\bbR$ 
such that~$(0,0,\ell)$ lies in~$C$.

To determine this, 
let $\xi_{m,n} = m(1-\iC\,x/y) + n(\iC/y)$.  
Thus, $\xi_{m,n}$ lies in~$\Lambda^{*}_{+}$ 
if and only if either $m>0$ or $m=0$ and $n>0$.  
Let $P$ be the plane in~$\bbR^{3}$ 
that contains the three (non-collinear!) points~$\hat\xi_{0,1}$,
$\hat\xi_{1,0}$, and~$\hat\xi_{1,1}$. 
 
 I claim that $\hat\xi_{m,n}$ lies strictly above $P$ for all 
$\xi_{m,n}\in \Lambda^{*}_{+}$ except when $(m,n)$ is one of $(0,1)$, $(1,0)$, or $(1,1)$.  (Here, `above' means on the opposite side of $P$ from $(0,0,0)\in\bbR^{3}$.) To see why, consider the determinant~$\Delta(m,n)$ of the matrix
$$
\begin{pmatrix}
\hat\xi_{0,1} & 1\\
\hat\xi_{1,0} & 1\\
\hat\xi_{1,1} & 1\\
\hat\xi_{m,n} & 1
\end{pmatrix}=
\begin{pmatrix} 
-1/y^{2} & 0 & 1/y^{2} & 1\\
1-x^{2}/y^{2} & -2x/y & 1+x^{2}/y^{2} & 1\\
1-(x{-}1)^{2}/y^{2} & -2(x{-}1)/y & 1+(x{-}1)^{2}/y^{2} & 1\\
m^{2}-(mx{-}n)^{2}/y^{2} & -2m(mx{-}n)/y & m^{2}+(mx{-}n)^{2}/y^{2} & 1
\end{pmatrix}.
$$
One finds that $\Delta(m,n) = (4/y^{3})\bigl(1-(m^{2}-mn+n^{2})\bigr)$.
The claim now follows by noting that $\Delta(0,0)>0$; that $\Delta(m,n) = 0$ 
if and only $(m,n)$ is one of $\pm(0,1)$, $\pm(1,0)$, or $\pm(1,1)$; 
and that $\Delta(m,n)<0$ otherwise.

Set $\xi_{1}=\xi_{0,1}$, $\xi_{2}=\xi_{1,0}$, and $\xi_{3}=\xi_{1,1}$.
It now follows that the triangle~$T$ 
with vertices $\hat\xi_{1}$, $\hat\xi_{2}$, $\hat\xi_{3}$ is a face of~$C$.  
Moreover, one finds that the ray~$\left\{(0,0,\ell)\ \vrule\ \ell\ge0\right\}$
meets~$T$ uniquely in the point
$$
r_{1}\,\hat\xi_{1}+r_{2}\,\hat\xi_{2}+r_{3}\,\hat\xi_{3}
=\bigl(0,0,S(x,y)\bigr)
$$
where $r_{1}\ge r_{2}\ge r_{3}\ge0$ satisfy $r_{1}+r_{2}+r_{3}=1$ 
and are given by
$$
r_{1} = \frac{y^{2}+x^{2}-x}{(y^{2}+x^{2}-x+1)}\,,\quad
r_{2} = \frac{1-x}{(y^{2}+x^{2}-x+1)}\,,\quad
r_{3} = \frac{x}{(y^{2}+x^{2}-x+1)},
$$
while 
$$
S(x,y) = \frac{2}{(y^{2}+x^{2}-x+1)}.
$$

Thus, it follows that
$$
V(\phi) \ge \frac{4\pi^{2}y}{(y^{2}+x^{2}-x+1)}\,,
$$
and that, moreover, equality holds only if $\phi_{\xi}=0$ 
for $\xi\not\in\left\{ \pm\xi_{1},\pm\xi_{2}, \pm\xi_{3}\right\}$
while $|\phi_{\xi_{i}}|^{2} = \tfrac12r_{i}$ for $i=1,2,3$.

To show that this lower bound is actually achieved by a $\d A$-balanced
weakly conformal immersion, one must show that it is possible to choose
the $\phi_{\xi_i}$ for $i=1,2,3$ so that the resulting $\phi$ is 
actually weakly conformal, i.e., so that $\phi_{z}\cdot \phi_{z}\equiv0$.
It is easy to see that necessary and sufficient conditions for this
are the relations
$$
\phi_{\xi_{i}}\cdot\phi_{\xi_{j}} = 0\ \text{for all $i$ and $j$}
$$
while 
$$
\phi_{\xi_{i}}\cdot\overline{\phi_{\xi_{j}}} = 0\ \text{for all $i\not=j$}.
$$
Of course, one must also require $|\phi_{\xi_{i}}|^{2} = \tfrac12r_{i}$ 
for $i=1,2,3$.  These equations can be satisfied 
(uniquely up to orthogonal transformation in $\bbE^{n+1}$ because all of the 
inner products have been specified) by the vectors
$$
\begin{aligned}
\phi_{\xi_{1}} &= \tfrac12\sqrt{r_{1}}(1,\iC,0,0,0,0),\\
\phi_{\xi_{2}} &= \tfrac12\sqrt{r_{2}}(0,0,1,\iC,0,0),\\
\phi_{\xi_{3}} &= \tfrac12\sqrt{r_{3}}(0,0,0,0,1,\iC).
\end{aligned}
$$
Using these vectors, one finds that the corresponding map $f_\tau:T\to S^{5}$
given by
$$
f_\tau(z) = \sum_{j=1}^{3} \left(\phi_{\xi_{j}}\,\E^{2\pi\iC\la\xi_{j},z\ra}
          +\overline{\phi_{\xi_{j}}}\,\E^{-2\pi\iC\la\xi_{j},z\ra}\right)
$$
satisfies
$$
\left|\d f_\tau\right|^{2} = \frac{4\pi^{2}}{(y^{2}+x^{2}-x+1)}\, |\d z|^{2},
$$
so this $f_\tau$ is a homothetic immersion of $T=\bbC/\Lambda$ 
with its flat conformal metric, and its induced area has the minimal value 
$$
V(f_\tau)  = \frac{4\pi^{2}y}{(y^{2}+x^{2}-x+1)} = V_{b}(5,T)\,.
$$
The uniqueness up to orthogonal transformation 
of this $V_{b}$-minimizer is now clear.
\end{proof}

\begin{remark}[Estimating the conformal volume]
The map $f_\tau$ constructed at the end of the proof 
of Theorem~\ref{thm: Vbcomputations}
is, of course, the map that Montiel and Ros denote as $\psi_{xy}$ 
(where $\tau = x + \iC\,y$) in \cite{Montiel-Ros}. 
An immediate consequence is their Corollary~6, i.e., that
$$
V_c\bigl(T,\bigl[|\mathrm{d}z|^2\bigl]\bigr) \ge V(f_\tau) 
= \frac{4\pi^{2}y}{(y^{2}+x^{2}-x+1)}\,.
$$
\end{remark}

\section{Conformal volume of the fundamental immersion}\label{sec: convol}
In this section, I study the conformal volume of the conformal immersion
of minimal $\d A$-balanced volume of the torus $T = \bbC/\Lambda$, where
$\Lambda\subset\bbC$ is the lattice generated by $1$ and $\tau = x +\iC\,y$
where $0\le x\le \tfrac12$ and $y\ge \sqrt{1-x^{2}}$.  The three
shortest elements in the dual lattice $\Lambda^{*}$ are $\xi_{1} = \iC/y$,
$\xi_{2} = 1 - \iC\,x/y$, and $\xi_{3} = \xi_{1}+\xi_{2} = 1 - \iC\,(x{-}1)/y$.

I will use $z = u+\iC\,v$ as the standard coordinate on $\bbC$, 
so that the differential $\d z$ is well-defined on $T$, 
with induced flat metric $|\d z|^{2} = \d z{\circ}\d\bar z$.

As was explained in~\S\ref{sec: balvol}, 
the fundamental immersion~$f_\tau:T\to S^{5}\subset \bbE^{6}$ is given by
\be
f_\tau(u+\iC\,v)
= \begin{pmatrix}\openup 10\jot
   \sqrt{r_{1}}\,\cos\bigl(2\pi\,v/y\bigr)\hfill\\
   \sqrt{r_{1}}\,\sin\bigl(2\pi\,v/y\bigr)\hfill\\
   \sqrt{r_{2}}\,\cos\bigl(2\pi\,(u - v\,x/y)\bigr)\hfill\\
   \sqrt{r_{2}}\,\sin\bigl(2\pi\,(u - v\,x/y)\bigr)\hfill\\
   \sqrt{r_{3}}\,\cos\bigl(2\pi\,(u - v\,(x{-}1)/y)\bigr)\\
   \sqrt{r_{3}}\,\sin\bigl(2\pi\,(u - v\,(x{-}1)/y)\bigr)
  \end{pmatrix}\,,
\ee
where
\be
r_{1} = \frac{y^{2}+x^{2}-x}{(y^{2}{+}x^{2}{-}x{+}1)}\,,\quad
r_{2} = \frac{1-x}{(y^{2}{+}x^{2}{-}x{+}1)}\,,\quad
r_{3} = \frac{x}{(y^{2}{+}x^{2}{-}x{+}1)}.
\ee
These numbers satisfy $1>r_{1}\ge r_{2}\ge r_{3}\ge0$ and $r_{1}+r_{2}+r_{3}=1$.
It is also worth noting that $r_{2}+r_{3}\le \tfrac23$, 
with equality if and only if $(x,y) = \bigl(\tfrac12,\tfrac12\sqrt{3}\bigr)$, 
in which case $r_{i} = \tfrac13$ for all~$i$.

The map $f_\tau$ is a homothetic conformal immersion, satisfying
$$
|\d f_\tau|^{2} = \frac{4\pi^{2}}{(y^{2}+x^{2}-x+1)}\, |\d z|^{2},
$$
and so its induced area is 
$$
V(f_\tau)  = \frac{4\pi^{2}y}{(y^{2}+x^{2}-x+1)}\,.
$$
Note that, if $x=0$, then $r_{3}=0$ 
and the immersion $f_\tau$ has image in $S^{3}\subset\bbE^{4}$, 
otherwise, it is linearly full in $S^{5}$.

\begin{theorem}\label{thm: MRconj}
If $r_{1}\le \tfrac23$, then $V_{c}(f_\tau) = V(f_\tau)$.
In particular, in this range, one has the equality
$$
V_c\bigl(\bbC/\Lambda,\bigl[|\mathrm{d}z|^2\bigl]\bigr) = V(f_\tau) 
= \frac{4\pi^{2}y}{(y^{2}+x^{2}-x+1)}
$$
\end{theorem}

\begin{remark} In~\cite[p.~165]{Montiel-Ros}, 
Montiel and Ros had conjectured 
that this equality holds when $r_{1}<\tfrac23$. 
(Their condition $(x-\tfrac12)^{2}+y^{2}<\tfrac94$ 
is equivalent to $r_{1}<\tfrac23$.)
\end{remark}

\begin{proof}
For simplicity, set
$$
(s,t) = \bigl(2\pi\,(u - v\,x/y),\ 2\pi\,v/y\bigr)
$$
so that $0\le s,t\le 2\pi$ is a fundamental rectangle for $\Lambda$,
with $y\,\d s \w \d t = 4\pi^{2}\,\d u\w\d v$.  In terms of $s$ and $t$,
$$
f_\tau(s,t) = \begin{pmatrix}
  \sqrt{r_1}\,\cos t\\ \sqrt{r_1}\,\sin t\\
  \sqrt{r_2}\,\cos s\\ \sqrt{r_2}\,\sin s\\
  \sqrt{r_3}\,\cos (s{+}t)\\ \sqrt{r_3}\,\sin (s{+}t) 
\end{pmatrix}
=\begin{pmatrix} 
\sqrt{r_1}\,\E^{\iC t}\\
\sqrt{r_2}\,\E^{\iC s}\\
\sqrt{r_3}\,\E^{\iC(s+t)}\\
\end{pmatrix}.
$$
Let $a = (a^{1},a^{2},\ldots,a^{6})\in\bbE^{6}$ satisfy $a{\cdot}a < 1$
and define
$$
V(a) = \frac{y}{(y^{2}{+}x^{2}{-}x{+}1)}
\int_{0}^{2\pi}\int_{0}^{2\pi}  
\frac{(1-a{\cdot}a)}{(1-a{\cdot}f_\tau(s,t))^{2}} \,\d s\ \d t.
$$
Then 
$$
V_{c}(f_\tau) = \sup\left\{\ V(a)\ \mid\ a{\cdot}a <1\ \right\}.
$$
Because $f_\tau$ is equivariant with respect to a $2$-torus of rotations 
in~$\SO(6)$, one can apply a rotation from this torus
to reduce to the cases in which $a_1,a_3\ge0$ and $a_2=a_4=0$,
while $a_5 = \cos\theta\,\sqrt{{a_5}^2+{a_6}^2}$ 
and $a_6 = -\sin\theta\,\sqrt{{a_5}^2+{a_6}^2}$ for some angle~$\theta$. 
Setting
$$
\lambda = a_1\sqrt{r_1}\,,\qquad \mu = a_3\sqrt{r_2}\,,\qquad
\nu = \sqrt{{a_5}^2+{a_6}^2}\,\sqrt{r_3}\,,
$$
one has $\lambda,\mu,\nu\ge0$ and $\lambda+\mu+\nu 
\le ({a_1}^2{+}{a_3}^2{+}{a_5}^2{+}{a_6}^2)^{1/2}(r_1{+}r_2{+}r_3)^{1/2}<1$.  

If $r_3\not=0$,%
\footnote{ When $r_3=0$, one can simply set $\nu=0$ and omit the term $\nu^2/r_3$.  
I will let the reader check that the following argument goes through 
\emph{mutatis mutandis} in the simpler case in which~$r_3=0$.  
Anyway, Montiel and Ros had already verified the conjecture when $r_3=0$.}
set
\be\label{eq: V_theta def}
\begin{aligned}
V_\theta(\lambda,\mu,\nu) &= \frac{y}{(y^{2}{+}x^{2}{-}x{+}1)}
\left(1-\frac{\lambda^2}{r_1}-\frac{\mu^2}{r_2}-\frac{\nu^2}{r_3}\right)\\
&\qquad \times \int_0^{2\pi}\int_0^{2\pi}\frac{\d s\,\d t}
{\bigl(1{-}\lambda\,\cos t{-}\mu\,\cos s{-}\nu\,\cos(s{+}t{+}\theta)\bigr)^2}
\end{aligned}
\ee
Let 
$$
D_+ = \left\{ (\lambda,\mu,\nu)\ \vrule\ \ \lambda,\mu,\nu\ge0,\quad 
 \frac{\lambda^2}{r_1} {+}\frac{\mu^2}{r_2}{+}\frac{\nu^2}{r_3} < 1\ \right\}
 \subset\bbR^3
$$
be the first octant of the ellipsoid. Then
$$
V_c(f_\tau) = \sup\left\{\ V_\theta(\lambda,\mu,\nu)
   \ \vrule\ \ \theta\in\bbR,\quad (\lambda,\mu,\nu)\in D_+\  \right\}.
$$
Considering the series expansion of the integral 
near $(\lambda,\mu,\nu)=(0,0,0)$, namely
$$
V_\theta(\lambda,\mu,\nu) = 
\frac{4\pi^2 y}{(y^{2}{+}x^{2}{-}x{+}1)}
\left(
\textstyle{1-\frac{\lambda^2}{r_1}-\frac{\mu^2}{r_2}-\frac{\nu^2}{r_3}}
\right)
\left(
\textstyle{1 + \frac32(\lambda^2+\mu^2+\nu^2) +6\lambda\mu\nu\cos\theta + \cdots}
\right),
$$
one sees that $(0,0,0)$ can be a local maximum of $V_\theta(\lambda,\mu,\nu)$
only if $r_i\le \frac23$ for $i=1,2,3$.  Since $0\le r_3\le r_2\le r_1$, 
the only condition this imposes is $r_1\le \frac23$.%
\footnote{Meanwhile, if $r_1>\frac23$, 
then $(\lambda,\mu,\nu)=(0,0,0)$ 
is \emph{not} a local maximum of $V_\theta(\lambda,\mu,\nu)$.}

Now, performing the $s$-integration%
\footnote{This uses the formula (valid for $a>|b|$)
$$\int_0^{2\pi} \frac{\d s}{(a - b\,\cos s)^2} 
= \frac{2\pi a}{(a^2-b^2)^{3/2}}.$$
}
in~\eqref{eq: V_theta def} yields
$$
\begin{aligned}
V_\theta(\lambda,\mu,\nu) &= \frac{2\pi y}{(y^{2}{+}x^{2}{-}x{+}1)}
\left(1-\frac{\lambda^2}{r_1}-\frac{\mu^2}{r_2}-\frac{\nu^2}{r_3}\right)\\
&\qquad \times
\int_0^{2\pi}\frac{(1-\lambda\,\cos t)\,\d t}
{\bigl((1{-}\lambda\,\cos t)^2
     {-}(\mu^2\,+2\mu\nu\,\cos(t{+}\theta)+\nu^2)\bigr)^{3/2}}\,.
\end{aligned}
$$
Since $\mu^2\,+2\mu\nu\,\cos(t{+}\theta)+\nu^2 \le (\mu{+}\nu)^2$, 
one has the inequality
\be\label{eq: lastintegral}
\begin{aligned}
V_\theta(\lambda,\mu,\nu) &\le \frac{2\pi y}{(y^{2}{+}x^{2}{-}x{+}1)}
\left(1-\frac{\lambda^2}{r_1}-\frac{\mu^2}{r_2}-\frac{\nu^2}{r_3}\right)\\
&\qquad\qquad \times
\int_0^{2\pi}\frac{(1-\lambda\,\cos t)\,\d t}
{\bigl((1{-}\lambda\,\cos t)^2 - (\mu{+}\nu)^2\bigr)^{3/2}}\,.
\end{aligned}
\ee

The remaining integral in~\eqref{eq: lastintegral} is not elementary, 
but it can be computed in terms of Legendre's complete elliptic integral
of the second kind~$E$,
which is the function defined as follows for $k\in[-1,1]$
$$
E(k) = \int_0^{\pi/2} \sqrt{1-k^2\,\sin^2\theta}\>\d\theta
\simeq\frac\pi2\left(1-\tfrac14 k^2 - \tfrac3{64}k^4 - \cdots\right)\,. 
$$
Because the coefficients of the convergent series expansion
are obviously negative for all positive $k$-degrees,
$E$ satisfies
\be\label{eq: Ebound}
E(k) \le \tfrac\pi2 (1-\tfrac14k^2),\qquad\text{$|k|\le 1$.}
\ee
Now, when $a,b\ge0$ and $a{+}b<1$, one has the formula%
\footnote{This integral formula is not entirely standard,
so I have included a derivation in Appendix~\ref{app: Ellipticintegral}.}
\be\label{eq: Ellipticintegral}
\int_0^{2\pi}\frac{(1-a\,\cos t)\,\d t}
{\bigl((1{-}a\,\cos t)^2 - b^2\bigr)^{3/2}}
= \frac{4\,E\left(\sqrt{{4ab}/{(1-(a{-}b)^2)}}\right)}
   {(1-(a{-}b)^2)^{1/2}(1-(a{+}b)^2)}\,.
\ee
Thus, \eqref{eq: lastintegral} becomes
$$
\begin{aligned}
V_\theta(\lambda,\mu,\nu) &\le \frac{8\pi y}{(y^{2}{+}x^{2}{-}x{+}1)}
\left(1-\frac{\lambda^2}{r_1}-\frac{\mu^2}{r_2}-\frac{\nu^2}{r_3}\right)\\
&\qquad \times
\frac{E\left(
  \sqrt{\mathstrut4\lambda(\mu{+}\nu)/(1-(\lambda{-}\mu{-}\nu)^2)}\right)}
{(1-(\lambda{-}\mu{-}\nu)^2)^{1/2}(1-(\lambda{+}\mu{+}\nu)^2)},
\end{aligned}
$$
which, using the inequality~\eqref{eq: Ebound}, implies
\be\label{eq: algebrabound}
\begin{aligned}
V_\theta(\lambda,\mu,\nu) &\le \frac{4\pi^2 y}{(y^{2}{+}x^{2}{-}x{+}1)}
\left(1-\frac{\lambda^2}{r_1}-\frac{\mu^2}{r_2}-\frac{\nu^2}{r_3}\right)\\
&\qquad \times
\frac{1-(\lambda{+}\mu{+}\nu)^2 + 3\lambda(\mu+\nu)}
{(1-(\lambda{-}\mu{-}\nu)^2)^{3/2}(1-(\lambda{+}\mu{+}\nu)^2)}.
\end{aligned}
\ee

Now, the right hand side of~\eqref{eq: algebrabound} will have its maximum
in~$D_+$ at $(\lambda,\mu,\nu)=(0,0,0)$ if $R\ge0$  in~$D_+$, where
$$
\begin{aligned}
R(\lambda,\mu,\nu) &= 
\bigl(1-(\lambda{+}\mu{+}\nu)^2\bigr)
\bigl(1-(\lambda{-}\mu{-}\nu)^2\bigr)^{3/2}\\
&\qquad 
- \left(1-\frac{\lambda^2}{r_1}-\frac{\mu^2}{r_2}-\frac{\nu^2}{r_3}\right)
   \bigl(1-(\lambda{+}\mu{+}\nu)^2 + 3\lambda(\mu+\nu)\bigr).\\
\end{aligned}
$$
Since $(1-c)^{3/2}\ge(1-\tfrac32c)$ when $0\le c\le 1$, 
one has $R\ge Q$  on $D_+$, where
$$
\begin{aligned}
Q(\lambda,\mu,\nu) &= 
\bigl(1-(\lambda{+}\mu{+}\nu)^2\bigr)
\bigl(1-\tfrac32(\lambda{-}\mu{-}\nu)^2\bigr)\\
&\qquad 
- \left(1-\frac{\lambda^2}{r_1}-\frac{\mu^2}{r_2}-\frac{\nu^2}{r_3}\right)
   \bigl(1-(\lambda{+}\mu{+}\nu)^2 + 3\lambda(\mu+\nu)\bigr)\\
&= Q_2(\lambda,\mu,\nu) + Q_4(\lambda,\mu,\nu),\\
&= \frac{\lambda^2}{r_1}+\frac{\mu^2}{r_2}+\frac{\nu^2}{r_3}
   - \frac32\bigl(\lambda^2+(\mu+\nu)^2\bigr) + Q_4(\lambda,\mu,\nu),
\end{aligned}
$$
and where $Q_j(\lambda,\mu,\nu)$ is homogeneous of degree~$j$.
Thus, it suffices to show that $Q\ge0$ on $D_+$.

Since $r_1\le\tfrac23$ and $r_2+r_3\le\tfrac23$,
the quadratic form $Q_2$ is nonnegative on~$\bbR^3$.  
Thus, when one considers, for $(\lambda,\mu,\nu)\in D_+$, 
the quartic polynomial in $t$ defined by
$$
q(t) = Q(t\lambda,t\mu,t\nu) 
= Q_2(\lambda,\mu,\nu)\,t^2 + Q_4(\lambda,\mu,\nu)\,t^4,
$$
one has $Q_2(\lambda,\mu,\nu)\ge0$.  
Meanwhile, assuming $(\lambda,\mu,\nu)\not=(0,0,0)$ 
and setting~$T = (\lambda{+}\mu{+}\nu)^{-1}>1$, 
one finds, directly from the defining formula of $Q$, that
$$
q(T) = T^4\left( 
\frac{\lambda^2}{r_1}+\frac{\mu^2}{r_2}+\frac{\nu^2}{r_3}
-(\lambda{+}\mu{+}\nu)^2\right)
\,\left(3\lambda(\mu{+}\nu)\right).
$$
Since $r_i<1$ while $r_1+r_2+r_3=1$, the quadratic form
$$
G(\lambda,\mu,\nu)
= \frac{\lambda^2}{r_1}+\frac{\mu^2}{r_2}+\frac{\nu^2}{r_3}
-(\lambda{+}\mu{+}\nu)^2
$$
is nonnegative. Thus, $q(T)\ge0$.   

Thus, writing $q(t)= t^2(a+bt^2)$, one has 
$$a = Q_2(\lambda,\mu,\nu)\ge0
\qquad\text{and}\qquad 
a+b\,T^2 = q(T)\,T^{-2}\ge0.
$$  
Since $T>1$, it follows that 
$$
Q(\lambda,\mu,\nu) = q(1) = a+b \ge 0.
$$

Thus, $R\ge Q\ge0$ on~$D_+$, implying that $V_\theta(\lambda,\mu,\nu)$
attains its supremum on~$D_+$ at~$(\lambda,\mu,\nu)=(0,0,0)$.  

Thus, when $r_1\le \tfrac23$, one has
$$
V_c(f_\tau) = \frac{4\pi^2 y}{(y^{2}{+}x^{2}{-}x{+}1)}
$$
as was to be shown.
\end{proof}

\begin{remark}[An upper bound on the conformal volume]\label{rem: confvol}
When $r_1>\tfrac23$, one can show (see Appendix~\ref{app: ConfVolCalc}) 
that the maximum value of $V_\theta(\lambda,\mu,\nu)$
occurs at $(\lambda,\mu,\nu) = \bigl(\sqrt{3r_1{-}2},0,0\bigr)$,
so
\be\label{eq: confvol f_tau}
V_c(f_\tau) = V_\theta\bigl(\sqrt{3r_1{-}2},0,0\bigr) 
= \frac{8\pi^2 y\sqrt{(y^2{+}x^2{-}x{+}1)}}{3\sqrt3\,(y^2{+}x^2{-}x)},
\ee
and this gives an upper bound for the conformal volume of $T=\mathbb{C}/\Lambda$.
For large $y$, this upper bound is about $20\%$ above the obvious
lower bound of $4\pi$.
\end{remark}

\begin{remark}[Relation with Willmore inequalities]
One of the original motivations for introducing the conformal volume
was to shed light on the Willmore Conjecture 
concerning the infimum of the integral
$$
W(\phi) = \int_T |H_\phi|^2\,\d A_\phi
$$
as $\phi$ ranges over the immersions of a torus~$T$ into $\bbE^n$.  
The conjecture (now proved~\cite{Marques-Neves}) was that $W(\phi)\ge 2\pi^2$.

The connection with conformal volume is that conformal volume
can be used to give a lower bound for the Willmore integral,
which is, itself, conformally invariant.  
In fact, as Li and Yau showed~\cite{Li-Yau}, 
when $\phi:\bbC/\Lambda\to \bbE^n$ is a conformal immersion,
one has
$$
W(\phi) \ge V_c\bigl(\bbC/\Lambda,\bigl[|\mathrm{d}z|^2\bigl]\bigr).
$$
Note that this alone, coupled with Theorem~\ref{thm: MRconj}, 
is enough to prove the Willmore conjecture for tori conformally
equivalent to~$\bbC/\Lambda$ 
when~$(y{-}1)^2+(x{-}\tfrac12)^2<\tfrac14$, 
as pointed out by Montiel and Ros~\cite[Corollary 7]{Montiel-Ros}.

Now, one could ask for a refined version of this lower bound,
one that fixes a conformal structure on the torus and allows~$\phi$ 
to range only over conformal immersions rather than all immersions.  

While a sharp form of such an inequality does not seem to be
within reach, by using arguments similar to those used to prove 
Theorem~\ref{thm: Vbcomputations}, one can prove, 
for any \emph{isometric} immersion~$\phi$ of $\bbC/\Lambda$
endowed with the metric $|\d z|^2$ into $\bbE^n$, 
that one has
$$
W(\phi) \ge \frac{\pi^2\bigl(y^4+(1{+}2x{-}2x^2)y^2-3x^2(1{-}x)^2\bigr)}{y^3}
$$
and that equality holds for such an isometric~$\phi$ 
if and only if $\phi=f_\tau$ up to translation and rotation.

Note that this lower bound is at least $2\pi^2$, 
with equality for $\tau=x + \iC\,y$ in the allowable hyperbolic triangle%
\footnote{I.e., $0\le x\le \tfrac12$ and $y\ge\sqrt{1{-}x^2}$.} 
if and only if $(x,y) = (0,1)$, i.e., for the square torus.
This lower bound also goes to infinity as $y\to\infty$.  
In particular, for large~$y$, it exceeds $8\pi$, 
which is the value of $W(\phi)$ when $\phi$ 
is a holomorphic branched double cover of the standard $2$-sphere.

For (unbranched) conformal embeddings~$\phi$, 
I am not aware of any violations of the above inequality.
While it is, perhaps, it is too much to hope that the above inequality 
holds for all conformal immersions (rather than just isometric ones), 
it is an interesting question as to what the lower bound
for conformal immersions might be.
\end{remark}

\appendix
\section{An elliptic integral}\label{app: Ellipticintegral}
Since the integral formula~\eqref{eq: Ellipticintegral} is not
entirely standard (in fact, none of the symbolic integration packages
to which I have access seem to be able to demonstrate it), I include
a proof here for the benefit of the reader.

Thus, assume that $a$ and $b$ are constants satisfying $a,b\ge0$ and $a{+}b<1$.
To evaluate the integral 
$$
I(a,b) = \int_0^{2\pi}\frac{(1-a\,\cos t)\,\d t}
{\bigl((1{-}a\,\cos t)^2 - b^2\bigr)^{3/2}}
= 2 \int_{0}^{\pi}\frac{(1-a\,\cos t)\,\d t}
{\bigl((1{-}a\,\cos t)^2 - b^2\bigr)^{3/2}},
$$
first, set
$$
\tau = \sqrt{\frac{1-b+a}{1-b-a}}
\qquad\text{and}\qquad
k = \sqrt{\frac{4ab}{1-(b-a)^2}}
$$
and note that $0\le k\le 1$. Now, make the substitution
$$
t = 2\tan^{-1}\left(\frac{\tan\theta}{\tau}\right),
$$
in the integral, noting that
$$
\cos t = \frac{\tau^2-\tan^2\theta}{\tau^2+\tan^2\theta}
\qquad\text{and}\qquad
\d t = \frac{2\tau\sec\theta\,\d\theta}{\tau^2+\tan^2\theta}.
$$
One then obtains
$$
\begin{aligned}
I(a,b) &= \frac{4\tau^2}{\bigl(1{-}(b{-}a)^2\bigr)^{3/2}(1{-}b{+}a)^2}
\\
&\quad \times
\int_0^{\pi/2}\frac{(1{-}a)(1{-}b{+}a)^2-2a\bigl(1{-}(b{-}a)^2\bigr)\sin^2\theta
+ 4a^2b\sin^4\theta}{(1-k^2\sin^2\theta)^{3/2}}\,\d\theta.
\end{aligned} 
$$
Now, let
$$
\omega_1 = \d\left(\frac{\sin\theta\cos\theta}{(1-k^2\sin^2\theta)^{1/2}}\right)
 = \frac{1-2\sin^2\theta+k^2\sin^4\theta}{(1-k^2\sin^2\theta)^{3/2}}\,\d\theta
$$
and
$$
\omega_2 = (1-k^2\sin^2\theta)^{1/2}\,\d\theta 
= \frac{1-2k^2\sin^2\theta+k^4\sin^4\theta}{(1-k^2\sin^2\theta)^{3/2}}\,\d\theta
$$
and note that, by the fundamental theorem of calculus and by definition, 
one has
$$
\int_0^{\pi/2}\omega_1 = 0
\qquad\text{and}\qquad
\int_0^{\pi/2}\omega_2 = E(k),
$$
where $E$ is Legendre's complete elliptic integral of the second kind.

Now, note the algebraic identity
$$
\frac{(1{-}a)(1{-}b{+}a)^2-2a\bigl(1{-}(b{-}a)^2\bigr)\sin^2\theta
+ 4a^2b\sin^4\theta}{(1-k^2\sin^2\theta)^{3/2}}\,\d\theta
= c_1\,\omega_1 + c_2\,\omega_2\,,
$$
where
$$
c_1 = \frac{a\bigl((1{-}b)^2{-}a^2\bigr)\bigl(1{-}(b{-}a)^2\bigr)}{(1{+}b)^2-a^2}
\qquad\text{and}\qquad
c_2 = \frac{\bigl(1{-}(b{-}a)^2\bigr)^2}{(1{+}b)^2-a^2}\,.
$$
Thus, 
$$
I(a,b) = \frac{4\tau^2}{\bigl(1{-}(b{-}a)^2\bigr)^{3/2}(1{-}b{+}a)^2}\ 
         \left(\frac{\bigl(1{-}(b{-}a)^2\bigr)^2}{(1{+}b)^2-a^2}\right)\ E(k),
$$
which, given the formulae for $\tau$ and $k$, simplifies to
$$
I(a,b) = \frac{4\,E\left(\sqrt{{4ab}/{(1-(a{-}b)^2)}}\right)}
           {(1-(a{-}b)^2)^{1/2}(1-(a{+}b)^2)},
$$
which is~\eqref{eq: Ellipticintegral}.

\section{Conformal Volume Computation}\label{app: ConfVolCalc}

Here is an argument to justify \eqref{eq: confvol f_tau}
by using~\eqref{eq: algebrabound}.

The first step is to use the algebraic equality
$$
\frac{\mu^2}{r_2} + \frac{\nu^2}{r_3}
 = \frac{(\mu+\nu)^2}{r_2{+}r_3} + \frac{(r_2\nu-r_3\mu)^2}{r_2r_3(r_2{+}r_3)}
$$
to deduce the inequality 
$$
\frac{\mu^2}{r_2} + \frac{\nu^2}{r_3}
 \ge \frac{(\mu+\nu)^2}{r_2{+}r_3} = \frac{(\mu+\nu)^2}{1{-}r_1}\,,
$$ 
which implies that the upper bound for $V_\theta(\lambda,\mu,\nu)$ 
provided by~\eqref{eq: algebrabound} can be \emph{weakened} to
\be\label{eq: weakalgebrabound}
\begin{aligned}
V_\theta(\lambda,\mu,\nu) 
&\le \frac{4\pi^2 y}{(y^{2}{+}x^{2}{-}x{+}1)}
\left(1-\frac{\lambda^2}{r_1}-\frac{(\mu{+}\nu)^2}{1{-}r_1}\right)\\
&\qquad \times
\frac{1-(\lambda{+}\mu{+}\nu)^2 + 3\lambda(\mu{+}\nu)}
{(1-(\lambda{-}\mu{-}\nu)^2)^{3/2}(1-(\lambda{+}\mu{+}\nu)^2)}.
\end{aligned}
\ee
Even though this is a weaker bound than~\eqref{eq: algebrabound}, it
has the advantages that $\mu$ and $\nu$ 
only occur as the combination $\mu{+}\nu$ and that only $r_1$ appears.
Thus, this right-hand expression is essentially a \emph{two}-variable
function depending on a \emph{single} parameter.

Thus, my goal is now to find the supremum of the function
$$
G_a(u,v) = \left(1-\frac{u^2}{a}-\frac{v^2}{1{-}a}\right)
\frac{1-(u{+}v)^2 + 3uv}
{(1-(u{-}v)^2)^{3/2}(1-(u{+}v)^2)}
\quad ( = G_{1-a}(v,u) )
$$
for $a$ satisfying $2/3 < a < 1$ in the region $Q_a$ 
in the $uv$-plane defined by the inequalities
$\frac{u^2}{a}+\frac{v^2}{1{-}a}<1$  and $u,v\ge0$.  

I claim that the supremum of $G_a$ in~$Q_a$ 
is attained only at the point $(u,v) = (3a{-}2,0)$. 
Moreover, one has
$$
G_a(3a{-}2,0) = \frac{2\sqrt3}{9a\sqrt{1{-}a}}\,.
$$
Once this has been shown, \eqref{eq: confvol f_tau} 
follows because then \eqref{eq: algebrabound} implies
$$
V_\theta(\lambda,\mu,\nu) \le V_0(\sqrt{3r_1{-}2},0,0),
$$ 
with equality if and only if $(\lambda,\mu,\nu) = (\sqrt{3r_1{-}2},0,0)$.

Unfortunately, $G_a$ does not extend continuously to the closure
of $Q_a$, i.e., the set $\overline{Q_a}$ defined by the inequalities 
$\frac{u^2}{a}+\frac{v^2}{1{-}a}\le1$ and $u,v\ge0$.  However, 
$G_a$ does extend (smoothly) to $\overline{Q_a}$ \emph{minus} 
the single point $(u,v) = (a,1{-}a)$.  

Now, it is not hard to show that, for every $\epsilon>0$,
there is a $\delta>0$ such that for $(u,v)\in Q_a$
satisfying $\bigl(u-a\bigr)^2+\bigl(v-(1-a)\bigr)^2<\delta$, 
the inequalities
$$
0 < G_a(u,v) \le \frac{3}{8\sqrt{a(1{-}a)}} + \epsilon
$$
hold. Since 
$$
\frac{3}{8\sqrt{a(1{-}a)}} < \frac{2\sqrt3}{9a\sqrt{1{-}a}}
$$
when $2/3 < a < 1$, it follows that the supremum of $G_a$ in $Q_a$
must occur in some compact set $R_{a,\delta}$ defined, for some $\delta>0$, 
by the inequalities $u\ge 0$, $v\ge 0$, 
$\bigl(u-a\bigr)^2+\bigl(v-(1-a)\bigr)^2\ge\delta$ and $\frac{u^2}{a}+\frac{v^2}{1{-}a}\le1$.

Obviously, $G_a$ vanishes on the points of the ellipse $\frac{u^2}{a}+\frac{v^2}{1{-}a}=1$ (except $(u,v)=(a,1{-}a)$, where it is not defined).
It is straightforward to show that on the boundary of $R_{a,\delta}$ 
where $v=0$, the unique maximum value of $G_a$ is at 
$(u,v) = (\sqrt{3a{-}2},0)$ and that on the boundary of $R_{a,\delta}$
where $u=0$, we have $G_a(u,0)\le 1$.

Thus, if the maximum value of $G_a$ on $R_{a,\delta}$ were greater than 
$G_a(\sqrt{3a{-}2},0)$, it would have to occur at a point $p_0=(u_0,v_0)$ 
in the interior of $R_{a,\delta}$, and, hence, $p_0$ would be a critical 
point of $G_a$.  

Thus, if one can show that there is no critical point of $G_a$ 
in the interior of $R_{a,\delta}$, it will follow that supremum of $G_a$ 
in $Q_a$ is attained at $(u,v) = (\sqrt{3a{-}2},0)$.  

To understand these critical points, consider the partials of $G_a$,
which take the form
$$
\frac{\partial G_a}{\partial u}
 = \frac{U(a,u,v)}{a(1{-}a)\bigl(1-(u{-}v)^2\bigr)^{5/2}\bigl(1-(u{+}v)^2\bigr)^2}
$$
and
$$
\frac{\partial G_a}{\partial v}
 = \frac{V(a,u,v)}{a(1{-}a)\bigl(1-(u{-}v)^2\bigr)^{5/2}\bigl(1-(u{+}v)^2\bigr)^2}\,,
$$
where $U(a,u,v)$ and $V(a,u,v)$ are polynomials in $a$, $u$ and $v$.
The critical points of $G_a$ in $Q_a$ occur only where $U(a,u,v)=V(a,u,v)=0$.

Computing\footnote{For this computation, I used MAPLE}
the Gröbner basis of the ideal generated 
by $U(a,u,v)$ and $V(a,u,v)$ with respect to the pure lexicographical order
$a > u > v$, 
one finds that the first element of this Gröbner basis is a polynomial
in $u$ and $v$ that factors as
$$
uv\bigl(1{-}(u{-}v)^2\bigr)^2\bigl(1{-}(u{+}v)^2\bigr)^2(1{-}u^2{+}uv{-}v^2)
\bigl(1{-}(u{-}v)^2(2{-}u^2{+}uv{-}v^2)\bigr).
$$

Now, it is easy to see that none of the factors of this polynomial 
vanish in the interior of $Q_a$.  Hence there are no critical points 
of~$G_a$ in the interior of $Q_a$.  Thus, there are no critical points 
of~$G_a$ in the interior of $R_{a,\delta}$, as was to be shown.

\bibliographystyle{hamsplain}

\providecommand{\bysame}{\leavevmode\hbox to3em{\hrulefill}\thinspace}

\end{document}